\documentclass{ifacconf}

\usepackage{graphicx}      
\usepackage{natbib}        
\usepackage{svg}           
\usepackage{adjustbox}     
\usepackage{amsmath}
\usepackage{amssymb}
\usepackage{enumerate}
\usepackage{bbm}

\begin{document}
\begin{frontmatter}

\title{Numerical Methods for Optimal Boundary Control of Advection-Diffusion-Reaction Systems} 


\author{Marcus Johan Schytt}, 
\author{John Bagterp J{\o}rgensen} 

\address{Department of Applied Mathematics and Computer Science, Technical University of Denmark, DK-2800 Kgs. Lyngby, Denmark (e-mail: \{mschytt, jbjo\}@dtu.dk)}

\begin{abstract}                
This paper considers the optimal boundary control of chemical systems described by advection-diffusion-reaction (ADR) equations. We use a discontinuous Galerkin finite element method (DG-FEM) for the spatial discretization of the governing partial differential equations, and the optimal control problem is directly discretized using multiple shooting. The temporal discretization and the corresponding sensitivity calculations are achieved by an explicit singly diagonally-implicit Runge Kutta (ESDIRK) method. ADR systems arise in process systems engineering and their operation can potentially be improved by nonlinear model predictive control (NMPC). We demonstrate a numerical approach for the solution to their optimal control problems (OCPs) in a chromatography case study. Preparative liquid chromatography is an important downstream process in biopharmaceutical manufacturing. We show that multi-step elution trajectories for batch processes can be optimized for economic objectives, providing superior performance compared to classical gradient elution trajectories.
\end{abstract}


\end{frontmatter}

\section{Introduction}
\label{sec:Introduction}
Advances in mechanistic modeling of labor-intensive, high-cost industrial processes have become a cornerstone for optimizing their operation. Among these models are systems of advection-diffusion-reaction (ADR) (or transport reaction) equations that describe the interplay between transport phenomena and reaction kinetics.
Given a mixture of $N_\mathcal{C}$ reactive components $\mathcal{C}$, the dynamics of their concentrations $c=[c_i]_{i \in \mathcal{C}}$ is described by the ADR system
\begin{align}\label{eq:ADRsystem}
    \partial_t c &= - \partial_z N + R(c),
\end{align}
in the interval domain $\Omega = [0,L]$. Here, the flux vector $N = N_a + N_d$ is split into its advective and diffusive parts
\begin{subequations}
\begin{align}\label{eq:fluxsplit}
N_a &= v c, \\
N_d &= -D \partial_z c,
\end{align}
\end{subequations}
with velocity $v$, and diffusion coefficient $D$. ADR systems are ubiquitous in chemical engineering, modeling a variety of unit operations such as catalytic fixed bed reactors and adsorption columns. For such applications, we consider the stoichiometric form of the reaction term
\begin{align}\label{eq:stoichiometricreaction}
    R(c) = \nu^T r(c),
\end{align}
with a constant stoichiometric matrix $\nu$ and reaction rates $r$. Altogether, the ADR system \eqref{eq:ADRsystem}-\eqref{eq:stoichiometricreaction} forms a system of semi-linear PDEs. Its formulation is complete with initial conditions and Danckwerts' inflow-outflow boundary conditions. The initial concentrations are given by a concentration profile
\begin{align}\label{eq:IC}
    c(0,z) = c_0(z), \quad \forall z \in \Omega,
\end{align}
and for all $t \geq 0$ the boundary fluxes satisfy
\begin{subequations}\label{eq:BCs}
\begin{align}
    \label{eq:Robin}
    N(t,0) &= v c_{\text{in}}(t),\\
    \label{Neumann}N(t,L) &= v c(t,L),
\end{align}
\end{subequations}
with $c_\text{in}$ denoting the vector of inlet concentrations. We permit systems with spatially stationary reactive components. Such components have vanishing concentration fluxes, eliminating their need for boundary conditions. However, the initial condition still applies, and the components enter the ADR system as ODEs.

In the process design of chemical reactors, the trajectory of the inlet composition can be taken as a design variable to optimize for economic objectives such as yield, purity, and productivity. This may be seen as a problem of optimal boundary control. We denote the control $u = c_\text{in}$ and consider the following optimal boundary control problem (OCP) in Bolza form
\begin{subequations}\label{eq:OCP}
\begin{gather}
    \begin{aligned}
        &\begin{split}
        \label{eq:OCPminpart}\min_{c,u} \phi = 
        &\int_0^T\int_\Omega \psi_1(c(t,z)) \,\mathrm{d}z \,\mathrm{d}t \\
        &\qquad+ \int_0^T \psi_2(c(t,L)) +\psi_3(u(t))\,\mathrm{d}t,
        \end{split}\\
        &\begin{split}
            \text{s.t. } \eqref{eq:ADRsystem}, \eqref{eq:IC}, \eqref{eq:BCs}, \text{ and}
        \end{split}
    \end{aligned}\\
    \label{eq:boxstate}
         \underline{c}(t,z) \leq c(t,z) \leq \overline{c}(t,z),\\
    \label{eq:boxinput}
        \underline{u}(t) \leq u(t) \leq \overline{u}(t),
\end{gather}
\end{subequations}
for all $(t,z) \in [0,T] \times \Omega$. Repeated solution of such OCPs lies at the heart of nonlinear model predictive control (NMPC). In the present context, MPC strategies may use online estimates of the upstream component composition \citep{Horsholt:etal:ECC:2019,Horsholt:etal:IFAC:2019,Hørsholt:2019} to compute optimal predictive controls for downstream processing. With this goal in mind, this paper studies the numerical solution of \eqref{eq:OCP}. In particular, we consider a compelling case study in preparative liquid chromatography, a downstream separation process in biopharmaceutical development. The literature on model-based optimization of chromatography processes has recently been reviewed by \cite{Kawajiri:2021}. A key challenge is the optimization of elution trajectories to ensure optimal separation of components. Conflicting economic objectives such as yield, purity, and productivity, complicate this task. The inherent complexity of the resulting multi-objective OCPs is addressed by simple (gradient) trajectory discretizations \citep{Karlsson:2004, Knutson:2015, Leweke:2018} and black-box optimization methods. The goal of this paper is to demonstrate the effectiveness of mathematical methods for optimal control to support model-based optimization of chemical processes.

\subsection{Paper overview}
The remainder of the paper is organized as follows: Section \ref{sec:SpatialDiscretization} describes the spatial discretization of the ADR system by the local discontinuous Galerkin (LDG) finite element method, a discontinuous Galerkin finite element method (DG-FEM). Section \ref{sec:MultipleShooting} describes the simultaneous discretization of the OCP using direct multiple shooting, an improvement upon sequential direct single shooting and a counterpart to orthogonal collocation. Section \ref{sec:Integration} describes an explicit singly diagonally-implicit Runge Kutta (ESDIRK) scheme, a temporal discretization strategy with good properties for efficient sensitivity calculation. Section \ref{sec:ChromatographyCaseStudy} five presents a chromatographic case study in protein purification, wherein the elution trajectory is optimized for different economic objectives. Conclusions are provided in Section \ref{sec:Conclusion}.

\section{Spatial semi-discretization}
\label{sec:SpatialDiscretization}
The method of lines (MOL) is commonly used for the numerical solution of ADR systems. It is a two-step approach that involves first a spatial semi-discretization followed by a temporal discretization (i.e., numerical integration) of the resulting ODE system. The spatial discretization is performed by several techniques. Low-order finite difference (FD) or finite volume (FV) methods are widely used in the chromatography literature \citep{Karlsson:2004, Degerman:2006, vonLieres:2010}. However they cause excessive numerical dissipation in the convection-dominated regime \citep{Enmark:2011}. In this paper, we apply a high-order DG-FEM. DG-FEMs have emerged in recent years as a preferred alternative for chromatography due to their accuracy and numerical efficiency \citep{Meyer:2020, Breuer:2023}.

We discretize the interval domain $\Omega$ into $N_\mathcal{E}$ uniform elements of length $h = L/N_\mathcal{E}$. We write $\Omega = \cup_{k=1}^{N_\mathcal{E}} \mathcal{E}_k$
with interval elements $\mathcal{E}_k = [z_k,z_{k+1}]$, on which the concentration of component $i \in \mathcal{C}$ is approximated by element-wise polynomials of degree at most $p \geq 1$. Using Lagrangian basis functions $\{\ell_j^k\}_{j=0}^p$ on element $\mathcal{E}_k$, we write
\begin{align}\label{eq:polyapprox}
c_i(t,z)|_{\mathcal{E}_k}  \approx c_{h,i}^k(t,z) = \sum_{j=0}^p c^k_{h,i,j}(t)\ell_{j}^k(z).
\end{align}
Note that the initial condition \eqref{eq:IC} is approximated similarly. The Lagrangian basis is characterized by a nodal set $\{\xi^k_j\}_{j=0}^p$ satisfying the interpolation property
\begin{align}
     c^k_{h,i}(t,\xi^k_j) = c^k_{h,i,j}(t), \quad \forall t \geq 0.
\end{align}
The nodal set satisfies a number of additional properties important for the numerical implementation. In particular, it is assumed to contain the element endpoints, which facilitates the prescription of boundary conditions and numerical fluxes. Details are discussed in Appendix \ref{sec:Appendix}.

\subsection{Elemental semi-discrete formulation}
Unlike continuous finite element approximations, each basis function is only supported on its corresponding element. As a result, approximations may be discontinuous across element boundaries. 
Importantly, this means that the approximations are element local, and indeed so is the resulting discretization of the ADR system. However, approximations are then multiply defined at element interfaces, which necessitates the introduction of numerical fluxes. To obtain a consistent discretization of the diffusive flux, we employ a mixed formulation. We introduce the diffusive flux as an auxiliary variable ${q} = {N}_d$, and we instead approximate the equivalent first-order system
\begin{subequations}\label{eq:firstordersystem}
\begin{align}
    \partial_t {c} &= - \partial_z ( {N}_a + {q}) + R(c),\\
    {q} &= N_d.
\end{align}
\end{subequations}
Based on the nodal DG-FEM methodology in \cite{Hes:07}, we insert the approximation \eqref{eq:polyapprox} in \eqref{eq:firstordersystem} and test the resulting system component- and element-wise. Applying the divergence theorem twice and inserting numerical fluxes, we obtain the following elemental semi-discrete formulation\footnote{With minor abuse of notation, the reaction term $R_i$ acts on the $N_\mathcal{C}$ coefficients associated to each of the $p+1$ basis functions on $\mathcal{E}_k$.}
\begin{subequations}\label{eq:implicitsd}
\begin{align}\label{eq:firstimplicit}
    M_h \partial_t \vec{c}_{h,i}^k = &- S_h (v \vec{c}_{h,i}^k) + \left[\left(v c_{h,i}^k - (v c_{h,i})^*\right)\boldsymbol{\ell}^k\right]_{z_k}^{z_{k+1}} \nonumber\\
    &- S_h \vec{q}_{h,i}^k + \left[\left(q_{h,i}^k - (q_{h,i})^*\right)\boldsymbol{\ell}^k\right]_{z_k}^{z_{k+1}}\nonumber\\
    &+M_h R_i(\vec{c}_h^k),
    \end{align}
    \begin{align}
    M_h \vec{q}_{h,i}^k =  - DS_h \vec{c}_{h,i}^k + D\left[\left(c_{h,i}^k - (c_{h,i})^*\right)\boldsymbol{\ell}^k\right]_{z_k}^{z_{k+1}}.
\end{align}
\end{subequations}
Here we define vectors of approximation coefficients
\begin{subequations}
\begin{align}
    \vec{c}_{h,i}^k &= [c_{h,i,j}^k]_{j = 0}^{p}, & \vec{q}_{h,i}^k = [q_{h,i,j}^k]_{j = 0}^{p},\\
    \vec{c}_h^k &= [\vec{c}_{h,i}^k]_{i \in \mathcal{C}}, & \vec{q}_h^k = [\vec{q}_{h,i}^k]_{i \in \mathcal{C}},
\end{align}
\end{subequations}
a vector of Lagrangian basis functions
\begin{align}
    \boldsymbol{\ell}^k=[\ell^k_j]_{j=0}^{p},
\end{align}
the elemental mass and stiffness matrices
\begin{subequations}\label{eq:elementalmats}
\begin{align}
    \label{eq:mass}M_{h,i,j} = \int_{[z_k,z_{k+1}]} \ell_i^k(z) \ell_j^k(z) \,\mathrm{d}z, \quad\\
    S_{h,i,j} = \int_{[z_k,z_{k+1}]} \ell_i^k(z) \partial_z\ell_j^k(z) \,\mathrm{d}z,
\end{align}
\end{subequations}
and advective, diffusive, and auxiliary numerical fluxes
\begin{align}\label{eq:numfluxes}
    (vc_{h,i})^*,\quad (q_{h,i})^*,\quad (c_{h,i})^*,
\end{align}
respectively. The numerical fluxes must be appropriately specified to enforce the boundary conditions~\eqref{eq:BCs} and to keep the discretization stable and conservative. They will be discussed in a moment. Remark that the elemental matrices \eqref{eq:elementalmats} are element invariant. Their efficient evaluation is shown in \cite{Hes:07}. In particular, inversion of the mass matrix \eqref{eq:mass} is computationally well-conditioned, and we can consider the now explicit semi-discrete formulation
\begin{subequations}\label{eq:explicitsd}
\begin{align}\label{eq:maineq}
    \partial_t \vec{c}_{h,i}^k = &M_h^{-1}\bigg(- S_h (v \vec{c}_{h,i}^k) + \left[\left(vc_{h,i}^k - (vc_{h,i})^*\right)\boldsymbol{\ell}^k\right]_{z_k}^{z_{k+1}} \nonumber\\
    &- S_h \vec{q}_{h,i}^k + \left[\left(q_{h,i}^k - (q_{h,i})^*\right)\boldsymbol{\ell}^k\right]_{z_k}^{z_{k+1}}\bigg)\nonumber\\
    &+R_i(\vec{c}_h^k),\\
    \label{eq:auxeq}\vec{q}_{h,i}^k = & M_h^{-1} \left(- DS_h\vec{c}_{h,i}^k + D\left[\left(c_{h,i}^k - (c_{h,i})^*\right)\boldsymbol{\ell}^k\right]_{z_k}^{z_{k+1}}\right).
\end{align}
\end{subequations}
This formulation may tempt the reader to immediately resolve the auxiliary variable \eqref{eq:auxeq} and include it in \eqref{eq:maineq}. However, for this to be possible, the auxiliary numerical flux must first be guaranteed to be independent of the auxiliary variable.

\subsection{Numerical flux scheme}
We employ the numerical fluxes originating from the LDG method introduced in \cite{Cockburn:98}. As usual, the advective flux is specified by upwinding across element boundaries. On the other hand, the diffusive and auxiliary fluxes are defined by opposing upwind and downwind directions, respectively. For $k = 2, \dots,N_\mathcal{E}$, we let
\begin{subequations}\label{eq:flux}
\begin{align}
(vc_{h,i})^*(t,z_{k}) &= v_ic_{i,h}^{k-1}(t,z_k),\\
(q_{h,i})^*(t,z_{k}) &= q_{h,i}^{k-1}(t,z_k),\\
(c_{h,i})^*(t,z_{k}) &= c_{h,i}^{k}(t,z_{k}).
\end{align}
\end{subequations}
At the domain boundaries, we define the numerical fluxes such that the boundary conditions are satisfied, in a weak sense, by the average boundary fluxes. We introduce numerical fluxes at the leftmost boundary $z=0$ and let
\begin{subequations}\label{eq:fluxBC}
\begin{align}
(v_ic_{h,i})^*(t,0) &= -c_{h,i}^{1}(t,0)+2v_ic_\text{in,i}(t),\\
(q_{h,i})^*(t,0) &= -q_{h,i}^{1}(t,0),\\
(c_{h,i})^*(t,0) &= c_{h,i}^{1}(t,0),\\
\intertext{and at the rightmost boundary $z=L$, we define fluxes}
(v_ic_{h,i})^*(t,L) &= c_{h,i}^{N_\mathcal{E}}(t,L),\\
(q_{h,i})^*(t,L) &= -q_{h,i}^{N_\mathcal{E}}(t,L),\\
(c_{h,i})^*(t,L) &= c_{h,i}^{N_\mathcal{E}}(t,L).
\end{align}
\end{subequations}
With this choice, the numerical fluxes \eqref{eq:flux}-\eqref{eq:fluxBC} are inserted into \eqref{eq:explicitsd} and we resolve the auxiliary variable which leaves the reduced formulation
\begin{align}\label{eq:elementalform}
    \nonumber\partial_t \vec{c}_{h,i}^k &= T^k_{h,-1}\vec{c}_{h,i}^{k-1} + T^k_{h,0}\vec{c}_{h,i}^k + T^k_{h,+1}\vec{c}_{h,i}^{k+1}\\
    &\quad+ c_\text{in,i}b_{h}^k + R_i(\Vec{c}_{h}^k),
\end{align}
for some appropriate matrices $T^k_{h,-1}, T^k_{h,0}, T^k_{h,+1}$ and vector $b^k_h$. This formulation exposes the three-element stencil of the LDG scheme, which has a smaller stencil compared to other locally resolvable flux schemes. We refer to \cite{Arn:02} for a review of DG-FEM schemes and their associated numerical fluxes.

\section{OCP discretization}
\label{sec:MultipleShooting}
We transform the OCP \eqref{eq:OCP} into a finite-dimensional nonlinear programming (NLP) problem. To this end, sequential methods (control discretization) such as direct single shooting have been used extensively for chromatographic optimization \citep{Degerman:2006, Ng:2012, Knutson:2015, Bock:2021, Cebulla:2023}. In contrast, simultaneous methods (state and control discretization) such as orthogonal collocation \citep{Biegler:1984} has only seen recent application in \cite{Holmqvist:2015, Holmqvist:2016}, despite offering increased robustness \citep{Biegler:2007}. Although collocation methods are effective, they face challenges in error control as they depend on fixed temporal discretizations. We apply direct multiple shooting \citep{Bock:1984} as an alternative approach that can leverage adaptive numerical integration to effectively manage the temporal discretization error. 

\subsection{Direct multiple shooting}
We consider $N_\mathcal{S}$ shooting intervals defined by equidistant sampling $t_k = kT_s$ with sampling interval $T_s = T / N_\mathcal{S}$. On each shooting interval $\mathcal{S}_k = [t_k, t_{k+1})$, we apply a zero-order hold (ZOH) discretization of the inlet trajectory $\Vec{u}=[u_k]_{k=0}^{N_\mathcal{S}-1}$ given
\begin{align*}
    c_\text{in}(t) = u(t) = u_k, \quad \forall t \in \mathcal{S}_k,
\end{align*}
such that the control function is piece-wise constant. The component concentrations are discretized using the method of Section \ref{sec:SpatialDiscretization}, and the resulting semi-discrete state $x=[\Vec{c}_h^k]_{k=1}^{N_\mathcal{E}}$ satisfies the ODE system
\begin{align}\label{eq:ODEsystem}
    \partial_t x = f(x,u), \quad x(0) = \Tilde{x}_0,
\end{align}
for right-hand side $f$ defined by \eqref{eq:elementalform} and initial condition $\Tilde{x}_0$ given by the finite element approximation of \eqref{eq:IC}. Temporal discretization results in the fully discrete state vector $\Vec{x} = [x_k]_{k=0}^{N_\mathcal{S}}$ defined by the initial condition
\begin{align}\label{eq:ICdiscrete}
    x_0 = \Tilde{x}_0,
\end{align}
and further by numerical integration
\begin{align}\label{eq:integration}
    x_{k+1} = F_k(x_k,u_k) \approx x_k + \int_{t_k}^{t_{k+1}} f(x(t),u_k) \, \mathrm{d}t.
\end{align}
This will be covered in Section \ref{sec:Integration}. The objective function in \eqref{eq:OCP} is discretized similarly with finite element approximations. We denote the spatially discrete integrand of the objective function
\begin{align}
\begin{aligned}
    \Psi(x(t), u(t)) \approx &\int_\Omega \psi_1(c(t,z) \,\mathrm{d}z\\ &\quad + \psi_2(c(t,L)) + \psi_3(u(t)), 
\end{aligned}
\end{align}
and after temporal discretization, we get the fully discrete objective function written as a sum over shooting intervals
\begin{align}\label{eq:}
    \Phi = \sum_{k=0}^{N_\mathcal{S}-1}\Phi_k(x_k,u_k) \approx \sum_{k=0}^{N_\mathcal{S}-1}\int_{t_k}^{t_k+1} \Psi(x(t),u(t))\,\mathrm{d}t.
\end{align}
The state and input constraints \eqref{eq:boxstate}-\eqref{eq:boxinput} are also discretized, and we get the multiple shooting NLP problem
\begin{subequations}\label{eq:MultipleShooting}
 \begin{gather}
    \begin{aligned}
        &\begin{split}
            \min_{\Vec{x},\Vec{u}} \Phi = \sum_{k=0}^{\mathcal{N}_s-1}\Phi_k(x_k,u_k),\\
        \end{split}\\
        &\begin{split}
            \text{s.t. } \eqref{eq:ICdiscrete},\eqref{eq:integration}, \text{ and}
        \end{split}
    \end{aligned}\\
    \underline{\Vec{x}} \leq \Vec{x} \leq \overline{\Vec{x}},\\
    \underline{\Vec{u}} \leq \Vec{u} \leq \overline{\Vec{u}}.
\end{gather}   
\end{subequations}

\section{Numerical integration and sensitivity computations}\label{sec:Integration}
The semi-discrete ADR system \eqref{eq:ODEsystem} exhibits stiffness, inherent in both the naturally stiff chemical reactions and the spatial discretization of diffusive transport with small diffusion coefficients. Efficient numerical integration of these systems therefore requires adaptive implicit solvers. Additionally, trajectory sensitivity information needs to be computed to effectively solve \eqref{eq:MultipleShooting}. To this end, software based on multi-step BDF methods, such as DASSL/DASPK \citep{Petzold:1982, Li:2000} and CVODES/IDAS \citep{Hindmarsh:2005}, has been used successfully in the literature. However, problems with discontinuities may be better handled by adaptive single-step Runge-Kutta methods. In particular, ESDIRK methods have shown good computational efficiency \citep{Kristensen:2004,Jorgensen:etal:ESDIRK:2018}, and have previously been applied in multiple shooting applications \citep{Capolei:2012}. 

\subsection{ESDIRK methods}
We apply the four-stage third-order ESDIRK method presented in \cite{Kristensen:2004} for the approximation of \eqref{eq:integration}. It is both stiffly accurate and L-stable. We summarize the scheme as follows. We define the intermediate state $x_k^n$ at $t = t_k^m$, which arises after $n$ integration steps in the shooting interval $\mathcal{S}_k$. We advance to the state $x_k^{n+1}$ at time $t_k^{m+1}=t_k^m+h$, $h$ being the step size provided by an adaptive step size controller, by solution of the equations
\begin{subequations}\label{eq:ESDIRK}
\begin{align}
    X_1 &= x_k^n,\\
    \Theta_i &= x_k^n + h\sum_{j=1}^{i-1} a_{ij} f(X_j,u_k)
    ,\\
    X_i &= \Theta_i + h\gamma f(X_i, u_k),\\
    x_{k}^{n+1} &= X_4,
\end{align}
\end{subequations}
for stages $i=2,3,4$. The explicit first step ensures high stage order, and the embedded local error estimate
\begin{align}\label{eq:Error}
    \text{err} = h\sum_{j=1}^{4}d_j f(X_j,u_k),
\end{align}
is of order four. The computational efficiency of {ESDIRK} methods lies in the solution of \eqref{eq:ESDIRK}, which can be performed in sequence. Since the semi-discrete system is large in dimension, this greatly reduces the numerical effort. The residual equation
\begin{align}
    0 = R_i(X_i) = X_i - \Theta_i - h\gamma f(X_i,u_k),  
\end{align}
is solved using an inexact Newton method, which assumes that the Jacobian $J_x f$ is constant across all stages
\begin{align}\label{eq:simplifying}
    J_x f(X_i,u_k) \approx J_x f_k^n = J_x f(x^{n}_k,u_k).
\end{align}
The Newton scheme is iterated to convergence via
\begin{subequations}\label{eq:Newtonscheme}
\begin{align}
    M^n_k\delta X_i^{m} = - R_i(X_i),\\
    X_i^{m+1} = X_i^{m} + \delta X_i^{m},
\end{align}
\end{subequations}
using the sparse LU factorization of the iteration matrix
\begin{align}\label{eq:iterationmatrix}
    M^n_k = I - h\gamma J_x f_k^n.
\end{align}
Consequently, both the number of Jacobian evaluations and LU factorizations are greatly reduced during integration, and accuracy is maintained by adaptively controlling convergence of \eqref{eq:Newtonscheme} and local error with \eqref{eq:Error}. We follow the implementation details outlined in \cite{Hairer:1996} for general implicit Runge Kutta methods.

\subsection{Sensitivity analysis}
Solving \eqref{eq:MultipleShooting} with gradient-based NLP algorithms requires both state and control sensitivities. The principle of internal numerical differentiation (IND) introduced by \cite{Bock:1981} presents a methodology for sensitivity generation which closely follows the numerical integration. Numerically exact sensitivities of \eqref{eq:integration} are obtained by differentiating the integration scheme \eqref{eq:ESDIRK} and its Newton scheme \eqref{eq:Newtonscheme} directly, a task that is difficult to implement for adaptive codes. A successful implementation (DAESOL) is described in \cite{Albersmeyer:2005}. Instead, we opt for a simpler yet computationally effective staggered direct approach.

When the $(n+1)$-st integration step has been accepted, the sensitivities are approximated by differentiation of \eqref{eq:ESDIRK} using the simplifying assumption \eqref{eq:simplifying}. Similar to before, the state sensitivities can be computed sequentially due to the formulation of ESDIRK methods. Sensitivities with respect to the previous step satisfy the equations
\begin{subequations}
\begin{align}
    J_{x_k^n}X_1 &= I,\\
    J_{x_k^n}\Theta_i &= I + h\sum_{j=1}^{i-1}a_{ij}(J_x f_k^n)(J_{x_k^n} X_j),\\
    \label{eq:specialboy} J_{x_k^n} X_i &= J_{x_k^n}\Theta_i + h\gamma (J_x f_k^n)(J_{x_k^n} X_i),\\
    A &= J_{x_k^n} X_4.
\end{align}
\end{subequations}
Note that \eqref{eq:specialboy} is equivalent to
\begin{align}
    M_k^n J_{x_k^n} X_i = I + h\sum_{j=1}^{i-1}a_{ij}(J_x f_k^n)(J_{x_k^n} X_j),
\end{align}
whose solution efficiently reuses the factorized iteration matrix. The sensitivity across the entire shooting interval 
\begin{align}
    A_k^{n} = J_{x_k} x_k^{n}
\end{align}
is advanced to the next step by use of the chain rule
\begin{align}
    A_k^{n+1} = A A_k^n.
\end{align}
The control sensitivities are derived analogously with the simplifying assumption
\begin{align}
    J_u f(X_i,u_k) \approx J_u f(X_k^n,u_k).
\end{align}

\section{Chromotography case study}
\label{sec:ChromatographyCaseStudy}
Chromatography unfolds in a column packed with porous particles forming the stationary phase. A component mixture, suspended in a liquid buffer (the mobile phase), courses through the column. Components adsorb to the stationary phase with varying affinities, resulting in distinct retention times forming the basis for separation. An eluent component, a mobile phase modifier (MPM), disrupts interactions between the mobile and stationary phases. By controlling the trajectory of the eluent, the separation is further enhanced.

We follow the case study presented in \cite{Karlsson:2004}, which considers the separation of proteins using ion-exchange chromatography. Immunoglobulin G (IgG) antibodies are separated from a mixture of IgG, bovine serum albumin (BSA), and myoglobin (Mb). Sodium chloride (NaCl) is introduced as an MPM and its trajectory is optimized for a variety of objectives. To this end, we first describe the model introduced in \cite{Karlsson:2004}.

\subsection{Chromatography model}
Mathematical models of chromatography are two-fold: a column model capturing transport through the interstitial volume, and a kinetic model describing adsorption to and desorption from the stationary phase. The literature on mathematical models for chromatography is extensive and is perhaps best summarized by \cite{Guiochon:2006}. In practice, one-dimensional reduced-order models considering only advection and diffusion are commonplace. 

We study a column model known as the equilibrium dispersive model. A column of length $L$, volume $V$, and porosity $\varepsilon$ is modeled axially by $\Omega = [0,L]$. The transport of protein mobile phase concentrations is governed by Darcy's law and is described by the ADR system \eqref{eq:ADRsystem} with boundary conditions \eqref{eq:BCs}, reactive components $\mathcal{C} = \mathcal{C}_q = \{\text{IgG}, \text{BSA}, \text{Mb}\}$, interstitial velocity $v = v_\text{in}$, and apparent diffusion coefficient $D = D_\text{app}$. The reaction term $R$ describes the adsorption of proteins on the stationary phase. We let $q$ denote the vector of stationary phase concentrations of bound proteins, and write
\begin{align}
    R(c) = - \phi \partial_t q,
\end{align}
where $\phi = (1-\varepsilon)/\varepsilon$ denotes the liquid volume fraction. Note that stationary phase concentrations may fit within the present framework of \eqref{eq:ADRsystem} by appending stationary components. We use \eqref{eq:stoichiometricreaction} with components $\mathcal{C}$ described by
\begin{align}
    c = [c_\text{NaCl}, c_\text{IgG}, c_\text{BSA}, c_\text{Mb}, q_\text{IgG}, q_\text{BSA}, q_\text{Mb}],
\end{align}
reaction rates $r= \partial_t q$, and stoichiometric matrix
\begin{align}
    \nu = \left[\begin{array}{ccc}
    0; -\phi I; & I
    \end{array}\right].
\end{align}
Table \ref{tb:Column} summarizes the column model parameters.

\begin{table}[tb]
\begin{center}
\caption{Column parameters}\label{tb:Column}
\begin{tabular}{ccccc}
$L$ {[}m{]}         & $V$ {[}m$^3${]}     & $\varepsilon$ {[}---{]} & $v$ {[}m/min{]} & $D$ {[}m$^2$/min{]} \\ \hline
$3.00\cdot 10^{-2}$ & $1.00\cdot 10^{-6}$ & 0.32                    & $3.00\cdot 10^{-2}$       & $5\cdot 10^{-6}$              
\end{tabular}
\end{center}
\end{table}

The kinetic model is given by a competitive multi-component Langmuir isotherm. Given $i \in \mathcal{C}_q$, we write
\begin{align}\label{eq:Langmuir}
    \partial_t q_i = k_{\text{ads},i}c_iq_{\text{max},i}\left(1-\sum_{j \in \mathcal{C}_q} \frac{q_j}{q_{\text{max},j}}\right) - k_{\text{des},i}q_i,
\end{align}
where $q_{\text{max},i}$ denotes its maximum stationary phase concentration, $k_{\text{ads},i}$ denotes its adsorption rate, and $k_{\text{des},i}$ denotes its desorption rate. The rates are parameterized in \cite{Karlsson:2004} to incorporate information about the MPM. The resulting rates are given by
\begin{subequations}\label{eq:MPM}
\begin{align}
    k_{\text{ads},i} &= k_{\text{ads}0,i} \exp({\gamma_i c_\text{NaCl}}),\\
    k_{\text{des},i} &= k_{\text{des}0,i} c_\text{NaCl}^{\beta_i}.\label{eq:MPMdes}
\end{align}
\end{subequations}
Here $k_{\text{ads}0,i}$ denotes the adsorption rate constant, $ k_{\text{des}0,i}$ denotes the desorption rate constant, $\gamma_i$ is a parameter which characterizes the protein hydrophobicity, and $\beta_i$ is a parameter describing its ion-exchange
characteristics. The eluent component is assumed to be inert, i.e. $\partial_t q_\text{NaCl} = 0$. We use the isotherm parameters shown in Table \ref{tb:Isotherm}, which were estimated in \cite{Karlsson:2004}.
\begin{table}[tb]
\begin{center}
\caption{Isotherm parameters}\label{tb:Isotherm}
\begin{tabular}{@{\hskip -0.025em}l@{\hskip -0.025em}cccc@{\hskip -0.025em}}
                                     & Unit  & IgG     & BSA     & Mb      \\ \hline
$q_{\text{max}, i}$ & [kmol/m$^3$]      & $5.40\cdot10^{-4}$ & $1.04\cdot10^{-3}$ & $7.50\cdot10^{-4}$ \\
$k_{\text{des},i}$ &[kmol/(m$^3$min)] & $3.00\cdot10^{3}$ & $3.00\cdot10^{3}$ & $3.00\cdot10^{3}$ \\
$k_{\text{ads},i}$ &[kmol/(m$^3$min)] & $2.31\cdot10^{6}$ & $1.76\cdot10^{5}$ & $5.00\cdot10^{6}$ \\
$\gamma_i$ &[m$^3$/kmol] & 0.00 & 0.00 & 0.00\\
$\beta_i$ &[---]                        & 1.12    & 3.20     & 0.61   
\end{tabular}
\end{center}
\end{table}

\subsection{Process description}
The purification step is performed as a batch process. Assuming the column is pre-equilibrated with a buffer volume, the process can be divided into three steps. First is the loading phase, wherein the mixture of components is injected into the column. Then the elution phase, which attempts to separate and elute a target component. Once the target component has passed through the column, the stripping phase is initiated and the column is purged.

\begin{table}[tb]
\begin{center}
\caption{Inlet concentrations}\label{tb:Inlets}
\begin{tabular}{ccccc}
           &Unit &  IgG      & BSA    & Mb    \\ \hline
 $c_\text{in,i}$ &[kmol/m$^3$] &$2.67\cdot 10^{-6}$     & $5.97\cdot 10^{-6}$       & $1.11\cdot 10^{-5}$ 
\end{tabular}
\end{center}
\end{table}

\begin{table}[tb]
\begin{center}
\caption{Process phases}\label{tb:Phases}
\begin{tabular}{cccccc}
  & Unit    & Load    & Elution  & Strip \\ \hline
Duration &[min]        & 8             & Variable & 6     \\
$c_\text{in,NaCl}$ & [kmol/m$^3$] & $9.00\cdot 10^{-3}$ & Variable & 1.00 
\end{tabular}
\end{center}
\end{table}

\subsection{Chromatographic optimization}
Separation quality is quantified by an objective function, of which there are several in chromatographic optimization. 
Process yield and purity become increasingly important when the target protein is very expensive. The yield of IgG is defined at the column outlet as
\begin{align}\label{eq:Yield}
    Y = \frac{\int_{t_1}^{t_2}c_{\text{IgG}} \,\mathrm{d}t}{t_{\text{load}}c_{\text{in},\text{IgG}}},
\end{align}
where $c_{\text{in},\text{IgG}}$ is given in Table \ref{tb:Inlets}, $t_{\text{load}}$ is the duration of the loading phase, and $t_1,t_2$ are cut times for fraction collection. The fraction collection duration is constrained by a given purity requirement. Outlet purity is defined as
\begin{align}
    \Pi = \frac{c_\text{IgG}}{\sum_{j \in \mathcal{C}_q} c_j},
\end{align}
and we seek $\Pi \geq 99\%$. However, the cost of materials (stationary phase) and the cost of time may also dominate the execution of the process. Under these circumstances, a productivity objective is more important. We define productivity considering only the cost of time
\begin{align}\label{eq:Productivity}
    P = \frac{Y}{t_\text{tot}},
\end{align}
where $t_\text{tot}$ is the total process duration. To keep our exposition short, we fix the loading and stripping times, keeping the elution time variable, and consider only open-loop optimal controls maximizing the IgG yield. Table \ref{tb:Phases} summarizes the process phases.

\begin{figure}[tb]
\begin{center}
\adjustbox{trim=0 0.55cm 0 0cm}{\includegraphics[width=8.4cm]{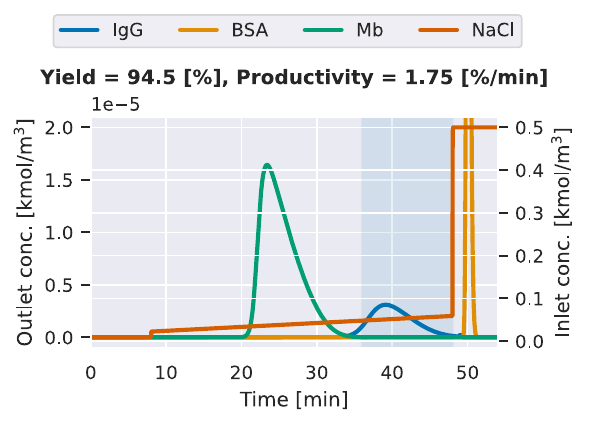}}    
\caption{Outlet profiles for optimal gradient elution with duration $T=40$ min.} 
\label{fig:subopt}
\end{center}
\end{figure}

The majority of literature on chromatographic optimization has shown linear (gradient) and single-step trajectories to be effective for separation. However, the constraint on the trajectory shape has the disadvantage of providing suboptimal solutions compared to multi-step or convex-concave trajectories (see \cite{Holmqvist:2016} and the references therein). Figure \ref{fig:subopt} shows an optimal gradient trajectory for the given process with fraction collection shaded in blue.

\begin{figure}[tb]
\begin{center}
\adjustbox{trim=0 0.55cm 0 0.0cm}{\includegraphics[width=8.4cm]{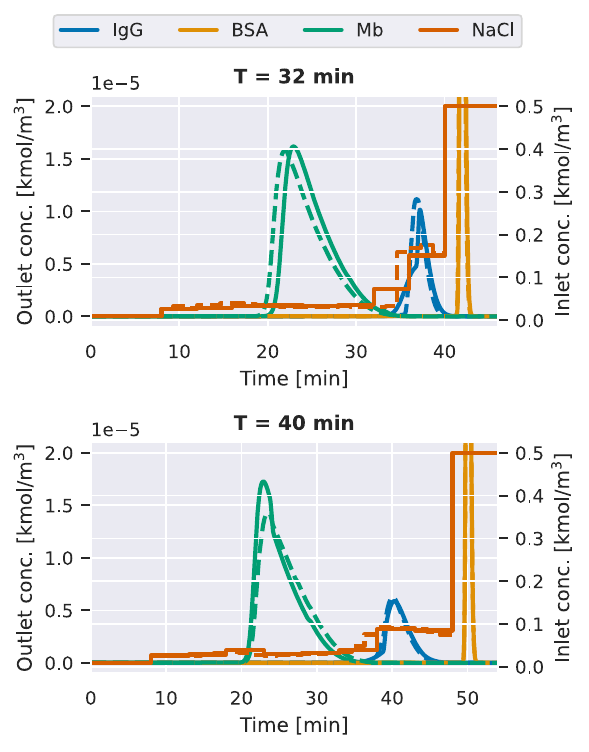}}  
\caption{Outlet profiles for optimal controls for $N_\mathcal{S} = 8$ (solid) and $N_\mathcal{S} = 24$ (dashed) shooting intervals.} 
\label{fig:optcont}
\end{center}
\end{figure}
\begin{figure}[tb]
\begin{center}
\adjustbox{trim=0 0.55cm 0 0cm}{\includegraphics[width=8.4cm]{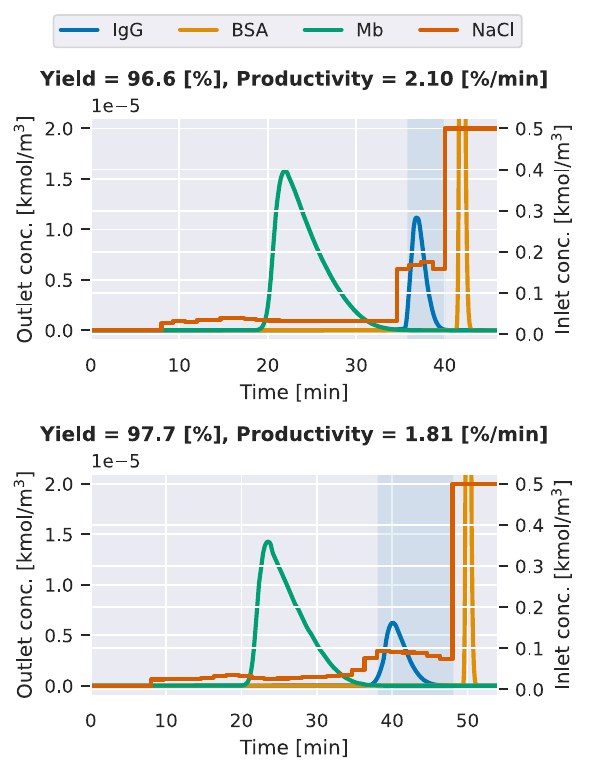}}  
\caption{Outlet profiles for optimal controls with duration $T=32$ min (top) and $T=40$ (bottom).} 
\label{fig:goodyield}
\end{center}
\end{figure}

\subsection{Numerical experiments}
We will improve upon the gradient elution trajectory using the numerical methods described in the previous sections. We discretize the column with DG-FEM, choosing third-order polynomials $p=3$, and $N_\mathcal{E} = 10$ elements as a good compromise between solution accuracy and problem dimension. ESDIRK is used for the numerical integration and sensitivity generation, and we choose a relative tolerance of $10^{-6}$ and an absolute tolerance of $10^{-8}$ to account for the problem scaling. Stricter tolerances are used to simulate the concentration profiles after optimization. We solve the OCP with multiple shooting, and the number of shooting intervals is chosen as $N_\mathcal{S} \in \{8, 24\}$, which results in NLP problems with 2528 and 7024 decision variables, respectively. We optimize with respect to yield \eqref{eq:Yield}, and report the resulting productivity with respect to elution durations $T \in \{32, 40\}$ min. To account for fraction collection, we implement smoothed exact cut times. We consider a sigmoid
\begin{align}
    \sigma(x) = 1/(1+\exp(-(x/\delta)),
\end{align}
with shape parameter $\delta>0$, and we define the objective function from \eqref{eq:OCP} with spatial term $\psi_1 = 0$ and outlet term
\begin{align}
    \psi_2 = -\frac{\sigma(\Pi - 0.99)c_\text{IgG}}{t_{\text{load}}c_{\text{in},\text{IgG}}},
\end{align}
such that we minimize the negative yield during fraction collection. We choose $\delta = 0.01$, and report the exact yield and productivity of the resulting optimal controls.

As is usual for NLP problems, it is important to consider problem scaling and the choice of good initial guesses. We solve the resulting NLP with IPOPT \citep{Wachter:2006}, which is equipped with automatic scaling options. To make the solver more robust, we include box-constraints that prevent states and controls from taking significant unphysical values which could make the integration fail. Our strategy for choosing good initial guesses is described as follows. The optimal gradient control is found and interpolated onto $N_\mathcal{S}$ shooting intervals. Using the initial loaded state, the elution trajectory is integrated with strict tolerances yielding initial shooting states. Since the problem suffers from suboptimal local minima, the interpolated control is randomly perturbed to obtain an ensemble of initial control guesses. We optimize these, seeking controls that satisfy a convergence tolerance of $\text{tol} = 10^{-8}$, and we choose the optimal control with the maximum yield. This procedure is then repeated for finer multiple shooting discretizations by interpolation of the newly found optimal control. The individual solution times can vary significantly, ranging in the order of several minutes, and we found that considering ensembles with hundreds of initial guesses was necessary.

The optimal multi-step trajectories are shown in Figure \ref{fig:optcont}. Compared to the optimal gradient trajectory for $T=40$ min, the IgG peaks are significantly higher and better separated from the Mb and BSA peaks. This can also be seen in the consequent yield and purity results shown in Figure \ref{fig:goodyield} for $N_\mathcal{S} = 24$, improving upon the optimal gradient trajectory. For $T=32$ min the optimal gradient trajectory degrades performance, recovering less than $86\%$ in the elution phase. In comparison, the optimal multi-step trajectory has near unchanged performance, showing a $-1.1\%$ yield loss and a $0.39${$\,\%/$min} productivity gain. This demonstrates the competitive nature of the yield and productivity objectives.

\section{Conclusion}
\label{sec:Conclusion}
This paper studies the numerical solution of OCPs governed by ADR systems. ADR systems are fundamental in the model-based optimization of chemical processes. We explore the spatial discretization of ADR systems using DG-FEM, a method which provides element-local formulations of the underlying PDE system. The OCPs are discretized by direct multiple shooting, and integrated in time using an ESDIRK method. ESDIRK methods provide efficient integration schemes from which approximate sensitivities may be generated without additional factorizations of the iteration matrix.

The use of DG-FEM and ESDIRK-based multiple shooting is exemplified in a chromatography case study in protein purification. A competitive Langmuir isotherm, specified with mobile phase modulating terms, is used to model the elution phase of a separation process. Optimal multi-step elution trajectories are found and show improved separation, providing better yield and productivity compared to conventional gradient trajectory strategies. However, to address the duration and difficulty of optimization, it is imperative to develop tailored implementations of the presented methods and utilize more classical global optimization techniques.


\bibliography{ifacconf}             

\begin{thebibliography}{36}
\providecommand{\natexlab}[1]{#1}
\providecommand{\url}[1]{\texttt{#1}}
\providecommand{\urlprefix}{URL }
\expandafter\ifx\csname urlstyle\endcsname\relax
  \providecommand{\doi}[1]{doi:\discretionary{}{}{}#1}\else
  \providecommand{\doi}{doi:\discretionary{}{}{}\begingroup \urlstyle{rm}\Url}\fi

\bibitem[{Albersmeyer(2005)}]{Albersmeyer:2005}
Albersmeyer, J. (2005).
\newblock \emph{Effiziente Ableitungserzeugung in einem adaptiven {BDF}-Verfahren}.
\newblock Master's thesis, Universit{\"a}t Heidelberg.

\bibitem[{Arnold et~al.(2002)Arnold, Brezzi, Cockburn, and Marini}]{Arn:02}
Arnold, D.N., Brezzi, F., Cockburn, B., and Marini, L.D. (2002).
\newblock Unified analysis of discontinuous galerkin methods for elliptic problems.
\newblock \emph{SIAM Journal on Numerical Analysis}, 39(5), 1749--1779.

\bibitem[{Biegler(1984)}]{Biegler:1984}
Biegler, L.T. (1984).
\newblock Solution of dynamic optimization problems by successive quadratic programming and orthogonal collocation.
\newblock \emph{Computers \& Chemical Engineering}, 8(3), 243--247.

\bibitem[{Biegler(2007)}]{Biegler:2007}
Biegler, L.T. (2007).
\newblock An overview of simultaneous strategies for dynamic optimization.
\newblock \emph{Chemical Engineering and Processing: Process Intensification}, 46(11), 1043--1053.

\bibitem[{Bock(1981)}]{Bock:1981}
Bock, H.G. (1981).
\newblock Numerical {{Treatment}} of {{Inverse Problems}} in {{Chemical Reaction Kinetics}}.
\newblock In K.H. Ebert, P.~Deuflhard, and W.~J{\"a}ger (eds.), \emph{Modelling of {{Chemical Reaction Systems}}}, Springer {{Series}} in {{Chemical Physics}}, 102--125. {Springer}, {Berlin, Heidelberg}.

\bibitem[{Bock and Plitt(1984)}]{Bock:1984}
Bock, H.G. and Plitt, K.J. (1984).
\newblock A {{Multiple Shooting Algorithm}} for {{Direct Solution}} of {{Optimal Control Problems}}.
\newblock \emph{IFAC Proceedings Volumes}, 17(2), 1603--1608.

\bibitem[{Bock et~al.(2021)Bock, Cebulla, Kirches, and Potschka}]{Bock:2021}
Bock, H.G., Cebulla, D.H., Kirches, C., and Potschka, A. (2021).
\newblock Mixed-integer optimal control for multimodal chromatography.
\newblock \emph{Computers \& Chemical Engineering}, 153, 107435.

\bibitem[{Breuer et~al.(2023)Breuer, Leweke, Schm{\"o}lder, Gassner, and Von~Lieres}]{Breuer:2023}
Breuer, J.M., Leweke, S., Schm{\"o}lder, J., Gassner, G., and Von~Lieres, E. (2023).
\newblock Spatial discontinuous {{Galerkin}} spectral element method for a family of chromatography models in {{CADET}}.
\newblock \emph{Computers \& Chemical Engineering}, 177, 108340.

\bibitem[{Capolei and J{\o}rgensen(2012)}]{Capolei:2012}
Capolei, A. and J{\o}rgensen, J.B. (2012).
\newblock Solution of {{Constrained Optimal Control Problems Using Multiple Shooting}} and {{ESDIRK Methods}}: {{American Control Conference}} ({{ACC}} 2012).
\newblock \emph{Proceedings of the 2012 American Control Conference}, 295--300.

\bibitem[{Cebulla et~al.(2023)Cebulla, Kirches, K{\"u}mmerer, and Potschka}]{Cebulla:2023}
Cebulla, D.H., Kirches, C., K{\"u}mmerer, N., and Potschka, A. (2023).
\newblock Model-based optimization of an ion exchange chromatography process for the separation of von {{Willebrand}} factor fragments and human serum albumin.
\newblock \emph{Proceedings in Applied Mathematics and Mechanics}, 00, e202300027.

\bibitem[{Cockburn and Shu(1998)}]{Cockburn:98}
Cockburn, B. and Shu, C.W. (1998).
\newblock The local discontinuous {G}alerkin method for time-dependent convection-diffusion systems.
\newblock \emph{SIAM Journal on Numerical Analysis}, 35(6), 2440--2463.

\bibitem[{Degerman et~al.(2006)Degerman, Jakobsson, and Nilsson}]{Degerman:2006}
Degerman, M., Jakobsson, N., and Nilsson, B. (2006).
\newblock Constrained optimization of a preparative ion-exchange step for antibody purification.
\newblock \emph{Journal of Chromatography A}, 1113(1-2), 92--100.

\bibitem[{Enmark et~al.(2011)Enmark, Arnell, Forss{\'e}n, Samuelsson, Kaczmarski, and Fornstedt}]{Enmark:2011}
Enmark, M., Arnell, R., Forss{\'e}n, P., Samuelsson, J., Kaczmarski, K., and Fornstedt, T. (2011).
\newblock A systematic investigation of algorithm impact in preparative chromatography with experimental verifications.
\newblock \emph{Journal of Chromatography A}, 1218(5), 662--672.

\bibitem[{Golub and Welsch(1969)}]{Gol:69}
Golub, G.H. and Welsch, J.H. (1969).
\newblock Calculation of gauss quadrature rules.
\newblock \emph{Mathematics of Computation}, 23(106), 221--230.

\bibitem[{Guiochon et~al.(2006)Guiochon, Felinger, Katti, and Shirazi}]{Guiochon:2006}
Guiochon, G.A., Felinger, A., Katti, A., and Shirazi, D.G. (2006).
\newblock \emph{Fundamentals of Preparative and Nonlinear Chromatography}.
\newblock {Elsevier, Amsterdam, Netherlands}.

\bibitem[{Hairer and Wanner(1996)}]{Hairer:1996}
Hairer, E. and Wanner, G. (1996).
\newblock \emph{Solving {{Ordinary Differential Equations II}}}, volume~14 of \emph{Springer {{Series}} in {{Computational Mathematics}}}.
\newblock {Springer}, {Berlin, Heidelberg}.

\bibitem[{Hesthaven(1997)}]{Hes:97}
Hesthaven, J. (1997).
\newblock From electrostatics to almost optimal nodal sets for polynomial interpolation in a simplex.
\newblock \emph{SIAM J. Numer. Anal.}, 35.

\bibitem[{Hesthaven and Warburton(2007)}]{Hes:07}
Hesthaven, J.S. and Warburton, T. (2007).
\newblock \emph{Nodal Discontinuous Galerkin Methods}.
\newblock Springer, New York.

\bibitem[{Hindmarsh et~al.(2005)Hindmarsh, Brown, Grant, Lee, Serban, Shumaker, and Woodward}]{Hindmarsh:2005}
Hindmarsh, A.C., Brown, P.N., Grant, K.E., Lee, S.L., Serban, R., Shumaker, D.E., and Woodward, C.S. (2005).
\newblock {SUNDIALS}: Suite of nonlinear and differential/algebraic equation solvers.
\newblock \emph{ACM Transactions on Mathematical Software (TOMS)}, 31(3), 363--396.

\bibitem[{Holmqvist and Magnusson(2016)}]{Holmqvist:2016}
Holmqvist, A. and Magnusson, F. (2016).
\newblock Open-loop optimal control of batch chromatographic separation processes using direct collocation.
\newblock \emph{Journal of Process Control}, 46, 55--74.

\bibitem[{Holmqvist et~al.(2015)Holmqvist, Andersson, Magnusson, and {\AA}kesson}]{Holmqvist:2015}
Holmqvist, A., Andersson, C., Magnusson, F., and {\AA}kesson, J. (2015).
\newblock Methods and {{Tools}} for {{Robust Optimal Control}} of {{Batch Chromatographic Separation Processes}}.
\newblock \emph{Processes}, 3(3), 568--606.

\bibitem[{Hørsholt et~al.(2019{\natexlab{a}})Hørsholt, Christiansen, Meyer, Huusom, and Jørgensen}]{Horsholt:etal:IFAC:2019}
Hørsholt, A., Christiansen, L.H., Meyer, K., Huusom, J.K., and Jørgensen, J.B. (2019{\natexlab{a}}).
\newblock A discontinuous-{G}alerkin finite-element method for simulation of packed bed chromatographic processes.
\newblock \emph{{IFAC-PapersOnline}}, 52(1), 346--351.

\bibitem[{Hørsholt et~al.(2019{\natexlab{b}})Hørsholt, Christiansen, Meyer, Huusom, and Jørgensen}]{Horsholt:etal:ECC:2019}
Hørsholt, A., Christiansen, L.H., Meyer, K., Huusom, J.K., and Jørgensen, J.B. (2019{\natexlab{b}}).
\newblock Spatial discretization and {K}alman filtering for ideal packed-bed chromatography.
\newblock In \emph{17th European Control Conference ({ECC})}, 2356--2361.

\bibitem[{Hørsholt et~al.(2019{\natexlab{c}})Hørsholt, Christiansen, Ritschel, Meyer, Huusom, and Jørgensen}]{Hørsholt:2019}
Hørsholt, A., Christiansen, L.H., Ritschel, T.K., Meyer, K., Huusom, J.K., and Jørgensen, J.B. (2019{\natexlab{c}}).
\newblock State and input estimation of nonlinear chromatographic processes.
\newblock In \emph{2019 IEEE Conference on Control Technology and Applications (CCTA)}, 1030--1035.

\bibitem[{Jørgensen et~al.(2018)Jørgensen, Kristensen, and Thomsen}]{Jorgensen:etal:ESDIRK:2018}
Jørgensen, J.B., Kristensen, M.R., and Thomsen, P.G. (2018).
\newblock A family of {ESDIRK} integration methods.
\newblock \emph{{arXiv:1803.01613}}.

\bibitem[{Karlsson et~al.(2004)Karlsson, Jakobsson, Axelsson, and Nilsson}]{Karlsson:2004}
Karlsson, D., Jakobsson, N., Axelsson, A., and Nilsson, B. (2004).
\newblock Model-based optimization of a preparative ion-exchange step for antibody purification.
\newblock \emph{Journal of Chromatography A}, 1055(1), 29--39.

\bibitem[{Kawajiri(2021)}]{Kawajiri:2021}
Kawajiri, Y. (2021).
\newblock Model-based optimization strategies for chromatographic processes: A review.
\newblock \emph{Adsorption}, 27(1), 1--26.

\bibitem[{Knutson et~al.(2015)Knutson, Holmqvist, and Nilsson}]{Knutson:2015}
Knutson, H.K., Holmqvist, A., and Nilsson, B. (2015).
\newblock Multi-objective optimization of chromatographic rare earth element separation.
\newblock \emph{Journal of Chromatography A}, 1416, 57--63.

\bibitem[{Kristensen et~al.(2004)Kristensen, J{\o}rgensen, Thomsen, and J{\o}rgensen}]{Kristensen:2004}
Kristensen, M.R., J{\o}rgensen, J.B., Thomsen, P.G., and J{\o}rgensen, S.B. (2004).
\newblock An {{ESDIRK}} method with sensitivity analysis capabilities.
\newblock \emph{Computers \& Chemical Engineering}, 28(12), 2695--2707.

\bibitem[{Leweke and {von Lieres}(2018)}]{Leweke:2018}
Leweke, S. and {von Lieres}, E. (2018).
\newblock Chromatography {{Analysis}} and {{Design Toolkit}} ({{CADET}}).
\newblock \emph{Computers \& Chemical Engineering}, 113, 274--294.

\bibitem[{Li and Petzold(2000)}]{Li:2000}
Li, S. and Petzold, L. (2000).
\newblock Software and algorithms for sensitivity analysis of large-scale differential algebraic systems.
\newblock \emph{Journal of Computational and Applied Mathematics}, 125(1), 131--145.

\bibitem[{Meyer et~al.(2020)Meyer, Leweke, Von~Lieres, Huusom, and Abildskov}]{Meyer:2020}
Meyer, K., Leweke, S., Von~Lieres, E., Huusom, J.K., and Abildskov, J. (2020).
\newblock {{ChromaTech}}: {{A}} discontinuous {{Galerkin}} spectral element simulator for preparative liquid chromatography.
\newblock \emph{Computers \& Chemical Engineering}, 141, 107012.

\bibitem[{Ng et~al.(2012)Ng, {Osuna-Sanchez}, Val{\'e}ry, S{\o}rensen, and Bracewell}]{Ng:2012}
Ng, C.K.S., {Osuna-Sanchez}, H., Val{\'e}ry, E., S{\o}rensen, E., and Bracewell, D.G. (2012).
\newblock Design of high productivity antibody capture by protein {{A}} chromatography using an integrated experimental and modeling approach.
\newblock \emph{Journal of Chromatography B}, 899, 116--126.

\bibitem[{Petzold(1982)}]{Petzold:1982}
Petzold, L.R. (1982).
\newblock Description of {{DASSL}}: A differential/algebraic system solver.
\newblock Technical Report SAND-82-8637; CONF-820810-21, {Sandia National Labs., Livermore, CA (USA)}.

\bibitem[{Von~Lieres and Andersson(2010)}]{vonLieres:2010}
Von~Lieres, E. and Andersson, J. (2010).
\newblock A fast and accurate solver for the general rate model of column liquid chromatography.
\newblock \emph{Computers \& Chemical Engineering}, 34(8), 1180--1191.

\bibitem[{W{\"a}chter and Biegler(2006)}]{Wachter:2006}
W{\"a}chter, A. and Biegler, L.T. (2006).
\newblock On the implementation of an interior-point filter line-search algorithm for large-scale nonlinear programming.
\newblock \emph{Mathematical Programming}, 106(1), 25--57.

\end{thebibliography}

\appendix
\section{Optimal nodal sets}\label{sec:Appendix}
The elemental nodal set on $\mathcal{E}_k$ arises for all $j=0,\dots,p,$ from an affine transformation
\begin{align}
    \xi_j^k = z_k + \frac{\xi_j + 1}{2},
\end{align}
of a reference nodal set $\{\xi_j\}_{j=0}^p$ on the reference element $\mathcal{E} = [-1,1]$. It is chosen with three properties in mind:
\begin{enumerate}[(i)]
  \item Minimize interpolation error.
  \item Associate with a quadrature rule.
  \item Have nodes at element boundaries.
\end{enumerate}
The importance of property (i) is obvious. Property (ii) is crucial, since the implementation of finite element methods requires accurate integration to avoid committing variational crimes. Finally, property (iii) simplifies the implementation of boundary conditions and numerical fluxes at element interfaces. We use the set of quadrature nodes associated to the Legendre-Gauss-Lobatto (LGL) rule as suggested in \cite{Hes:07}. It satisfies properties (ii)-(iii), and has near optimal interpolation error bounds as shown in \cite{Hes:97}. The LGL nodes consists of the $p+1$ roots of
\begin{align}\label{eq:findroots}
    \zeta(\xi)= (1-\xi^2)\partial_\xi P^{(0,0)}_{p}(\xi),
\end{align}
where $P^{(\alpha,\beta)}_n$ is the $n$-th degree Jacobi polynomial with weights $\alpha,\beta$.  Since the derivative of a Jacobi polynomial is another Jacobi polynomial of lower order, the roots of \eqref{eq:findroots} are found by first taking the boundary nodes as the endpoints $\xi = -1,1$, and then taking the $p-1$ interior nodes as the roots of
\begin{align}\label{eq:findrootsv2}
    \gamma(\xi) = \sqrt{p(p+1)}P^{(1,1)}_{p-1}(\xi).
\end{align}
The roots of \eqref{eq:findrootsv2} arise as the solution to an eigenvalue problem which can be solved numerically \citep{Gol:69}.

\end{document}